1

# On the Geometric Interpretation of the Complex Fourier Transforms of a Class of Exponential Functions

Jeremy Williams, London, England
E-mail: universaltutors@aol.com

Abstract: A class of complex Fourier transforms of exponential functions which have all their zeros on the real line is explored from a geometric perspective. These transforms belong to the Laguerre-Pólya class, and it is proved that all the zeros are simple.

Functions which map the complex plane to the complex plane provide an interesting challenge when it comes to visualisation. A full plot would require 4 dimensions, yet on paper we are obliged to draw in 2 dimensions. In this paper, we are going to explore ways of visualising a certain class of holomorphic entire functions defined for integer values of n ≥ 1, namely:

$$F_{2n}(z) = \int_{-\infty}^{\infty} e^{-t^{2n}} \cdot e^{izt} dt = R(\sigma,w) + i\, I(\sigma,w), \text{ where } z = w - i\sigma, \text{ and } \sigma, w, R, I \in \Re \qquad (1.0)$$

In 1923, Professor Pólya proved that all the zeros of these functions lie on the real line Im (z) = 0  (reference I).

**Case n = 1**    It is straightforward to show by contour integration that:

$$F_2(z) = \sqrt{\pi} \cdot e^{-z^2/4} = \sqrt{\pi} \cdot e^{(\sigma+iw)^2/4} = \sqrt{\pi} \cdot e^{(\sigma^2-w^2)/4} \cdot e^{(i\sigma w)/2}.$$

This is a Laguerre - Pólya function of order 2 in z. Whilst it is an elementary expression, we are going to examine some geometric aspects as they provide clues to the more involved cases which follow when n ≥ 2 .

The modulus of F squared, call it L (for length) squared: $L^2 = R^2 + I^2 = \pi \cdot e^{(\sigma^2-w^2)/2}$, so there are no zeros in this case. For constant $\sigma$ not equal to zero, and w = vt where v is a constant velocity and t represents time, a parametric plot of (R, I) yields a spiral orbit of rapidly diminishing radius:

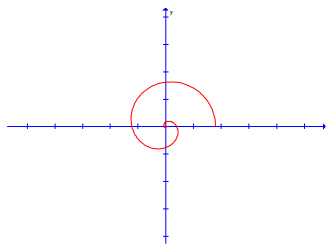

(Fig. 1, not to scale, $\sigma > 0$)

If a larger value of $\sigma$ is chosen, the orbit is enlarged.

The field lines R= 0 are given by w = π (1+2m)/ $\sigma$ and the field lines I = 0 are given by w= π (2m)/ $\sigma$, where m is any integer:

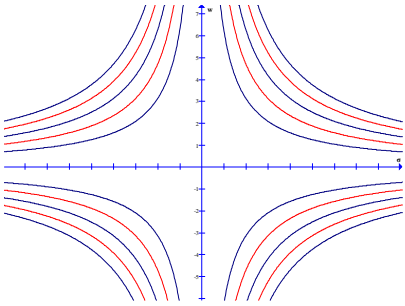

(Fig. 2)

We observe that the R= 0 field lines never meet each other. The same is true for the I = 0 field lines, and also no I = 0 field line ever meets an R= 0 field line (or else there would be a zero, which would establish a contradiction). A plot of $L^2$ against $\sigma$ for constant w rises steadily as $\sigma$ moves away from 0:

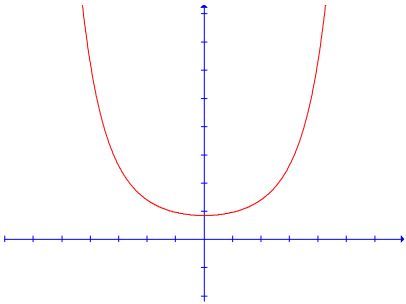

(Fig. 3)

Whereas a plot of $L^2$ against w for constant $\sigma$ falls steadily as w moves away from 0:

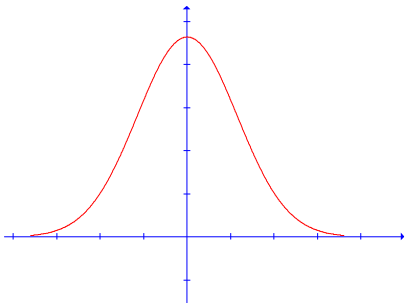

(Fig. 4)

## Case  n ≥ 2

In the same paper previously mentioned (I), Pólya also proved that in this case there are an infinite number of real zeros. In another paper (II) published 4 years later, Pólya proved that the order of such functions is at most $2n/(2n-1) < 2$.

By Hadamard's factorisation theorem (see for example III), such (entire) functions have the product representation:

$F_{2n}(z) = z^m e^{G(z)} \prod_{k=1}^{\infty} (1 - z/w_k)$, where G(z) is a polynomial of degree at most d, where d $\geq 0$ is an integer such that d $\leq v \leq$ d +1, and $v$ is the order of $F_{2n}(z)$. As $v < 2$, we see d < 2, and so the most general form of $F_{2n}(z)$ is given by $F_{2n}(z) = z^m e^{(a+b.z)} \cdot \prod_{k=1}^{\infty}(1-z/w_k) = c.z^m \cdot e^{bz} \prod_{k=1}^{\infty}(1-z/w_k)$    (1.1)

We now observe that z = 0 $\Rightarrow F_{2n}(z) = \int_{-\infty}^{\infty} e^{-t^{2n}} dt > 0$ and so m = 0 and c > 0. As the exponential kernel in the integral $\int_{-\infty}^{\infty} e^{-t^{2n}} \cdot e^{izt} dt$ is even in t, it follows that for every root there is another root of equal magnitude and opposite sign.

If we pair each root $w_k$ with its partner of opposite sign, then $F_{2n}(z) = c \cdot e^{bz} \prod_{r=1}^{\infty}(1 - z^2/\alpha_r^2)$ where $w_1 = \alpha_1$, $w_2 = -\alpha_1$, $w_3 = \alpha_2$, $w_4 = -\alpha_2$ etc, and the $\alpha_r$ are the positive roots. Owing to the kernel $e^{-t^{2n}}$ in the integral expression (equation 1.0) being even in t, it follows that $F_{2n}(z)$ is an even function in z, and so b = 0.

Accordingly, we see that: $F_{2n}(z) = c \prod_{k=1}^{\infty}(1 - z/w_k)$  (1.2)

and also that $F_{2n}(z) = c \cdot \prod_{r=1}^{\infty}(1 - z^2/\alpha_r^2)$  (1.3)

These functions belong to the so-called Laguerre-Pólya class, which we will denote by L-P.

Now let us consider $L^2 = F_{2n}(z) \cdot \overline{F_{2n}(z)} = c^2 \prod_{r=1}^{\infty}(1 - z^2/\alpha_r^2)(1 - (\overline{z})^2/\alpha_r^2)$

As z = w – i σ, $L^2 = c^2 \prod_{r=1}^{\infty}(1 - ((w - i\sigma)/\alpha_r)^2)(1 - ((w + i\sigma)/\alpha_r)^2)$

The rth term inside the product, call it $p_r$, = 1 - $2(w^2 - \sigma^2)/\alpha_r^2 + (w^2 + \sigma^2)^2/\alpha_r^4$, which can be written $(1 - (w^2 + \sigma^2)/\alpha_r^2)^2 + 4\sigma^2/\alpha_r^2$, which is always positive for positive σ. When σ = 0, $p_r$ reduces to $(1 - w^2/\alpha_r^2)^2$. Now $\frac{\partial p_r}{\partial \sigma}$ is $4\sigma/\alpha_r^2 + 4\sigma(\sigma^2 + w^2)/\alpha_r^4$ which is positive for positive σ. If σ = 0 and w is equal to one of the roots, then $L^2 = 0$ at that point. Since the zeros of holomorphic functions are isolated, as we increase σ above 0 (keeping w constant) there must immediately be an interval where $L^2 > 0$, and as each $p_r$ is strictly increasing for positive σ, so $L^2$ must increase as σ increases, and keep increasing. If on the other hand we start with σ = 0 and pick a value of w which is not a root, then $L^2$ is greater than 0 anyway, and as we increase σ above zero (with w constant), $L^2$ must increase as each $p_r$ is strictly increasing for positive σ.

These arguments show that the graph of $L^2$ against σ for constant w has a similar form to the one in Fig.3 previously, except that if we start with w = a root, the curve passes through the origin:

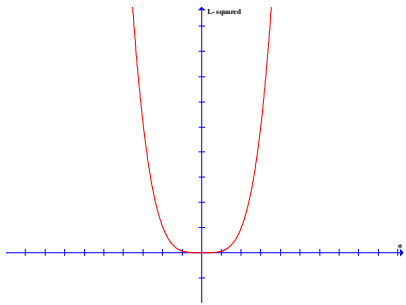
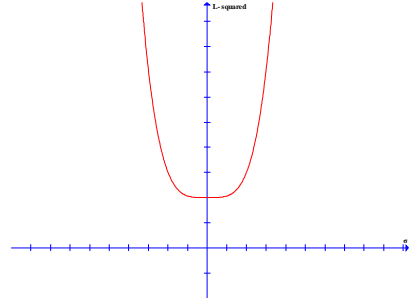

Fig. 5: $L^2$ against σ, constant w (a root )   Fig. 6: $L^2$ against σ, constant w (not a root)

This reveals one geometric reason why all the zeros of $F_{2n}(z)$ as defined by equation (1.0) lie on the real line: for any given constant w, the modulus of the function rises monotonically if we travel away from the real line $\sigma = 0$.

I will now explore $L^2$ from a different angle – by looking at the integral formula which defines $F_{2n}(z)$:

$$F_{2n}(z) = \int_{-\infty}^{\infty} e^{-t^{2n}} \cdot e^{izt} dt \quad \text{where } z = w - i\sigma$$

Now $L^2 = F_{2n}(z) \cdot \overline{F_{2n}(z)}$

$$= (\int_{-\infty}^{\infty} e^{(-x^{2n}+\sigma.x)}(\cos wx + i\sin wx)dx) \cdot (\int_{-\infty}^{\infty} e^{(-y^{2n}+\sigma.y)}(\cos wy - i\sin wy)dy)$$

$$= \iint_{\Re^2} e^{(-x^{2n}-y^{2n}+\sigma.(x+y))} \cdot ((\cos wx.\cos wy + \sin wx.\sin wy) + i(\sin wx.\cos wy - \cos wx.\sin wy)).dxdy$$

$$= \iint_{\Re^2} e^{(-x^{2n}-y^{2n}+\sigma.(x+y))} \cdot (\cos(wx-wy) + i\sin(wx-wy))dxdy$$

$$= \iint_{\Re^2} e^{(-x^{2n}-y^{2n}+\sigma.(x+y)+i.w(x-y))} dxdy \qquad \text{(manipulations justified by Fubini's Theorem).}$$

Now as $e^{\sigma(x+y)} = \sum_{m \geq 0} \sigma^m (x+y)^m / m!$ we can write $L^2$ as:

$$\sum_{m \geq 0} \frac{\sigma^m}{m!} \iint_{\Re^2} (x+y)^m \cdot e^{-(x^{2n}-y^{2n})+iw(x-y)} dxdy \quad \text{(justified by absolute uniform convergence)}$$

and substituting $t = x + y$, $X = x - y$ we see:

$$L^2 = \sum_{m \geq 0} \frac{\sigma^m}{m!} \iint_{\Re^2} t^m e^{-((t+X)/2)^{2n}+((t-X)/2)^{2n})+iwX} \left\| \frac{\partial(x,y)}{\partial(t,X)} \right\| dtdX$$

$$= \frac{1}{2} \sum_{m \geq 0} \frac{\sigma^m}{m!} \iint_{\Re^2} t^m e^{-((t+X)/2)^{2n}+((t-X)/2)^{2n})+iwX} dtdX$$

and as $e^{-(((t+X)/2)^{2n}+((t-X)/2)^{2n})}$ is even in t, for odd powers of n the integrals must be zero, so we need only consider the even powers and we can write:

$$L^2 = \frac{1}{2} \sum_{m \geq 0} \frac{\sigma^{2m}}{(2m)!} A_{2m,2n}(w) \qquad (1.4)$$

where $\quad A_{2m,2n}(w) = \iint_{\Re^2} t^{2m} e^{-(((t+X)/2)^{2n}+((t-X)/2)^{2n})+iwX} dtdX \qquad (1.5)$

## Lemma 1

$A_{2m,2n}(w) \geq 0 \qquad \forall\ m \geq 0,\ n \geq 2$

## Proof of Lemma 1

$$A_{2m,2n}(w) = \frac{1}{i^{2m}} \cdot \frac{\partial^{2m}}{\partial u^{2m}} \left[ \iint_{\Re^2} e^{-(((t+X)/2)^{2n} + ((t-X)/2)^{2n}) + i(wX+ut)} . dt dX \right]_{u=0}$$

and as the imaginary part of the exponent in the integrand can be expressed as

(w+u)(X+t)/2 +(u-w)(t-X)/2, if we now reverse the substitution t = x + y, X = x-y we obtain:

$$A_{2m,2n}(w) = (-1)^m \cdot \frac{\partial^{2m}}{\partial u^{2m}} \left[ \iint_{\Re^2} e^{-(x^{2n}+y^{2n}) + i((w+u)x + (u-w)y)} . 2 dx dy \right]_{u=0}$$

$$= (-1)^m \cdot \frac{\partial^{2m}}{\partial u^{2m}} \left[ \int_{-\infty}^{\infty} e^{-x^{2n}+ix(w+u)} dx \cdot \int_{-\infty}^{\infty} e^{-y^{2n}+iy(u-w)} dy . 2 \right]_{u=0} \qquad \text{by Fubini's theorem}$$

$$= (-1)^m . 2 . \frac{\partial^{2m}}{\partial u^{2m}} \left[ F_{2n}(u+w) . F_{2n}(u-w) \right]_{u=0} \qquad (1.6)$$

Now let us pick one particular value of n. Recall equation (1.3), namely that $F_{2n}(z) = c . \prod_{r=1}^{\infty}(1 - z^2/\alpha_r^2)$. Consider the partial product $P_N(z) = c . \prod_{r=1}^{N}(1 - z^2/\alpha_r^2)$ and suppose that for one particular value of N, K say, $T_{K,m}(w) = (-1)^m 2 . \frac{\partial^{2m}}{\partial u^{2m}} \left[ P_K(u+w) . P_K(u-w) \right]_{u=0} \geq 0$ $\forall\ m \geq 0,\ \forall\ w \in \Re$.

Now let us consider $T_{K+1,m}(w) = (-1)^m 2 . \frac{\partial^{2m}}{\partial u^{2m}} \left[ P_{K+1}(u+w) . P_{K+1}(u-w) \right]_{u=0}$

$$= (-1)^m 2 . \frac{\partial^{2m}}{\partial u^{2m}} \left[ (1-(u+w)^2/\alpha_{K+1}^2)(1-(u-w)^2/\alpha_{K+1}^2) . P_K(u+w) . P_K(u-w) \right]_{u=0}$$

$$= (-1)^m 2 . \frac{\partial^{2m}}{\partial u^{2m}} \left[ (1 - 2(u^2+w^2)/\alpha_{K+1}^2 + (u^4+w^4-2u^2w^2)/\alpha_{K+1}^4) . P_K(u+w) . P_K(u-w) \right]_{u=0}$$

For brevity of writing, call $1 - 2(u^2+w^2)/\alpha_{K+1}^2 + (u^4+w^4-2u^2w^2)/\alpha_{K+1}^4 = B$ and

$P_K(u+w) . P_K(u-w) = D$

If m = 0 then $T_{K+1,0}(w) = 2\ P_{K+1}^2(w) \geq 0$

If m = 1 then $T_{K+1,1}(w) = -2 \left(B_{uu} D + 2 B_u D_u + B D_{uu}\right)_{u=0}$ where the subscripts denote partial differentiation with respect to (wrt) u. As B is an even function in u, the middle term is zero.

Now $B_{uu} = -4/\alpha_{K+1}^2 + (12 u^2 - 4 w^2)/\alpha_{K+1}^4$, so we obtain:

$$T_{K+1,1}(w) = 4 \left(1/\alpha_{K+1}^2 + w^2/\alpha_{K+1}^4\right) T_{K,0}(w) + (1 - w^2/\alpha_{K+1}^2)^2 \cdot T_{K,1}(w) \geq 0$$

If m $\geq$ 2 then $(-1)^m T_{K+1,m}(w)/2 = \dfrac{\partial^{2m}}{\partial u^{2m}}(B.D)_{u=0} =$

$$\left(B \frac{\partial^{2m} D}{\partial u^{2m}} + {}^{2m}C_1 B_u \cdot \frac{\partial^{2m-1} D}{\partial u^{2m-1}} + {}^{2m}C_2 \cdot B_{uu} \frac{\partial^{2m-2} D}{\partial u^{2m-2}} + {}^{2m}C_3 \cdot B_{uuu} \frac{\partial^{2m-3} D}{\partial u^{2m-3}} + {}^{2m}C_4 B_{uuuu} \cdot \frac{\partial^{2m-2} D}{\partial u^{2m-2}}\right)_{u=0}$$

As before, the odd terms are zero as B is an even function in u, and also partial derivatives of B w.r.t u of degree more than 4 are identically zero.

Then as additionally $B_{uuuu} = 24/\alpha_{K+1}^4$, we see $T_{K+1,m}(w) =$

$$\left(B \frac{\partial^{2m}}{\partial u^{2m}}[(-1)^m 2D] + {}^{2m}C_2 \cdot B_{uu} \frac{\partial^{2m-2}}{\partial u^{2m-2}}[(-1)^m 2D] + {}^{2m}C_4 B_{uuuu} \cdot \frac{\partial^{2m-4}}{\partial u^{2m-4}}[(-1)^m 2D]\right)_{u=0}$$

$$= \left(B \frac{\partial^{2m}}{\partial u^{2m}}[(-1)^m 2D] - {}^{2m}C_2 \cdot B_{uu} \frac{\partial^{2m-2}}{\partial u^{2m-2}}[(-1)^{m-1} 2D] + {}^{2m}C_4 B_{uuuu} \cdot \frac{\partial^{2m-4}}{\partial u^{2m-4}}[(-1)^{m-2} 2D]\right)_{u=0}$$

$$= (1 - \frac{w^2}{\alpha_{K+1}^2})^2 \cdot T_{K,m}(w) + {}^{2m}C_2 . 4(1/\alpha_{K+1}^2 + w^2/\alpha_{K+1}^4) \cdot T_{K,m-1}(w) + {}^{2m}C_4 . 24/\alpha_{K+1}^4 \cdot T_{K,m-2}(w)$$

which is $\geq 0$.

What has been established is the inductive step, namely that $T_{K+1,m}(w) \geq 0$.
If we now consider the case K = 1 it is easy to see that $T_{1,m}(w) \geq 0 \ \forall \text{ m} \geq 0$:

$T_{1,0} = 2 P_1^2(w) \geq 0$, $T_{1,1}(w) = 8 (1/\alpha_1^2 + w^2/\alpha_1^4)$, $T_{1,2}(w) = 48/\alpha_1^4$, $T_{1,M}(w)$ is 0 for M $\geq$ 3.

Therefore by induction $T_{N,m}(w) \geq 0$ for N as large as we like, $\forall$ m $\geq$ 0.

Therefore taking the limit as N $\to \infty$ shows the expression for $A_{2m,2n}(w)$ in equation (1.6)

is $\geq$ 0. This argument holds for all integer values of n from 2 upwards.

QED.

Looking again at formula for $L^2$ in (1.4), the graphs in Fig.5 and Fig.6 are now very clear, along with the result that there are no zeros of $F_{2n}(z)$ off the line $\sigma = 0$.

As the functions $A_{2m,2n}(w)$ are Fourier transforms and so tend to zero as $w \to \infty$, we can also see that the plot of $L^2$ against w for constant σ must also have the form of the curve in Fig.4 (perhaps with some minor undulations, as $A_{0,2n}(w) = F_{2n}^2(w) \geq 0$ has undulations because of the zeros $w_k$):

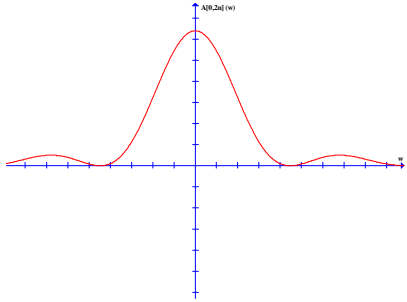

(Fig. 7)   Plot of $A_{0,2n}(w)$ against w.

Now we will consider the parametric plot of (R,I) when σ is constant and not equal to zero, and w = v t where v is a constant velocity.

The angular momentum $\underline{J} = m \ \underline{r} \wedge \underline{\dot{r}} = m \begin{vmatrix} e_1 & e_2 & e_3 \\ R & I & 0 \\ \dot{R} & \dot{I} & 0 \end{vmatrix}$   where the dots denote differentiation with respect to time, and the vectors in the top row are a left-handed set of orthonormal vectors, with the first being in the positive R direction and the second in the positive I direction.

so $\underline{J} = m \ (R.\dfrac{dI}{dw}.\dfrac{dw}{dt} - I.\dfrac{dR}{dw}.\dfrac{dw}{dt})\underline{e_3}$                 (1.7)

As we are holding σ fixed, we see that $\underline{J} = mv(R.\dfrac{\partial I}{\partial w} - I.\dfrac{\partial R}{\partial w})\underline{e_3}$

which by the Cauchy-Riemann equations can be expressed as :

$mv\left(R.\dfrac{\partial R}{\partial \sigma} + I.\dfrac{\partial I}{\partial \sigma}\right)\underline{e_3}$   ie $mv.\dfrac{\partial}{\partial \sigma}[R^2 + I^2]\underline{e_3}$          (1.8)

By equation (1.4), we can immediately see that for positive σ the angular momentum is also strictly positive, and as the $A_{2m,2n}(w)$ are Fourier transforms and tend to zero as $w \to \infty$, <u>the locus of (R,I) must follow a spiral orbit where the radius is never zero for any finite w</u>, as in Fig.1 at the beginning of this paper. (However, the orbit will swing in and out more than in the simple case in the first part of the paper.)

It is worth observing that, given any fixed value of n,  it is not possible as m varies for all the $A_{2m,2n}(w)$ to be zero at the same time for any particular value of w, w* say, as by equation (1.4) this would lead to a continuous line of zeros of $F_{2n}(z)$ in the complex plane along w = w*, which is impossible as holomorphic functions only have isolated zeros.

We are now going to explore the plot of the field lines R = 0 and I = 0:

Trivially, I = 0 when w = 0, because the integral in equation (1.0) is then real, and I = 0 when σ = 0, as the integral kernel is then even in t.

Now consider the spiral orbit discussed on the last page. This means that as we travel up the line σ = constant at speed v, we are crossing the field lines R = 0, then I = 0, R = 0 etc.

The field lines with equation R = 0 traverse the zeros on the w – axis.

If we consider a field line R = 0, then along this line we have:

$$R_\sigma + R_w \frac{dw}{d\sigma} = 0, \text{ or alternatively } R_\sigma \frac{d\sigma}{dw} + R_w = 0 \qquad (1.9)$$

where the suffixes denote the first partial differentials.

Now L-P is closed under differentiation, so it follows that all the zeros of the first derivative $F_{2n}^{(1)}(z)$ must also lie on the line σ = 0 too.

This means $R_w$ and $R_\sigma$ cannot both be zero at the same time off the line σ = 0, or else by the Cauchy-Riemann equations $F_{2n}^{(1)}(z)$ would have a zero, which would establish a contradiction.

Therefore, as we move along R= 0, unless $\frac{dw}{d\sigma}$ or $\frac{d\sigma}{dw}$ is zero, by (1.9), neither $R_w$ nor $R_\sigma$ can be zero (if one were zero the other would also have to be zero), which means that they must keep the same sign, and therefore so must $\frac{dw}{d\sigma}$.

All these factors considered, we might expect the R = 0 field line plot to look like that in Fig.2, except with the lines bent so that they cross the w–axis, intercepting the zeros, with the I = 0 field lines in between the R = 0 field lines when $\sigma$ is not equal to 0.

By way of example, consider the case n = 2, where:

$$\text{Re}[F_4(z)] = \int_{-\infty}^{\infty} e^{-t^4+\sigma \cdot t}\cdot\cos(wt)dt = 0 \qquad (2.0)$$

Consider positive values of $\sigma$. Letting t = T + k, and choosing k so that the coefficient of T in the exponent of the integrand vanishes leads us to set k = $(\sigma/4)^{(1/3)}$ and leads to:

$$\int_\Re e^{-(T^4+4T^3(\sigma/4)^{1/3}+6T^2(\sigma/4)^{2/3})}\cdot(\cos(wT)\cos(w(\sigma/4)^{1/3}) - \sin(wT)\sin(w(\sigma/4)^{1/3}))dT = 0$$

Replacing T with $u.(4/\sigma)^{1/3}$ then leads to:

$$\int_\Re e^{-(u^4.(4/\sigma)^{4/3}+4u^3(4/\sigma)^{2/3}+6u^2)}\cdot(\cos wu(4/\sigma)^{1/3}\cos w(\sigma/4)^{1/3} - \sin wu(4/\sigma)^{1/3}\sin w(\sigma/4)^{1/3})du$$
$$= 0.$$

Rearrangement then yields:

$$\cos(w(\sigma/4)^{1/3})\int_{\Re} e^{-(u^4 \cdot (4/\sigma)^{4/3} + 4u^3(4/\sigma)^{2/3} + 6u^2)} \cdot \cos(wu(4/\sigma)^{1/3})du =$$

$$\sin(w(\sigma/4)^{1/3})\int_{\Re} e^{-(u^4 \cdot (4/\sigma)^{4/3} + 4u^3(4/\sigma)^{2/3} + 6u^2)} \sin(wu(4/\sigma)^{1/3})du \qquad (2.1)$$

We can see that in the asymptotic case as $\sigma \to \infty$, the exponent in the integrand tends to $-6u^2$, so there is a family of asymptotic solutions given by w = $(4/\sigma)^{1/3} \cdot \frac{\pi}{2} \cdot (2m+1)$, where m is any integer. In the asymptotic limit equation (2.1) then reduces to:

0 = limit as $\sigma \to \infty$ [$(-1)^m \int_{\Re} e^{-(6u^2)} \sin(u(4/\sigma)^{2/3} \cdot \frac{\pi}{2} \cdot (2m+1))du$ ]; the sine function is odd and the exponential function is even, hence the integral is zero.

A plot of the R = 0 lines is shown below in Fig.8 (solid lines), along with part of the curves w = $(4/\sigma)^{1/3} \cdot \frac{\pi}{2} \cdot (2m+1)$ (dashed lines). As we can see, the asymptotic behaviour quickly becomes apparent for relatively small values of $\sigma$. The horizontal axis is the $\sigma$ axis, and the vertical axis is the w axis:

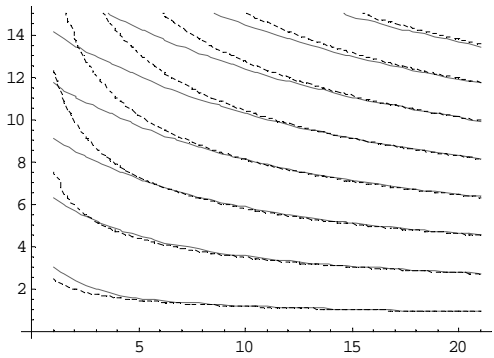

(Fig. 8)

Finally, a plot of R = 0 field lines (black) and I = 0 field lines (green) is shown in Fig. 9 below.

The zeros of $F_4(z)$ occur only where the R = 0 field lines cross the w axis (which is itself one of the I = 0 field lines), so all the zeros lie on the line $\sigma = 0$. As before, the horizontal axis is the $\sigma$ axis, and the vertical axis is the w axis:

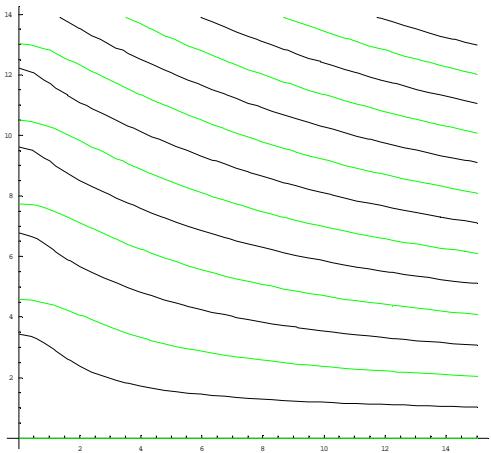

(Fig. 9)

A similar analysis can be carried out in the cases n ≥ 3, but the asymptotic convergence of the R = 0 field lines to the curves w = $(2n/\sigma)^{1/(2n-1)} \cdot \frac{\pi}{2}.(1+2m)$ as $\sigma \to \infty$ is slower.

We can see that there is another geometric reason for the location of the zeros of $F_{2n}(z)$: the only place where the R = 0 and I = 0 field lines intersect is on the w axis, ie where $\sigma = 0$.

## Simplicity of zeros

It will now be proved that all the zeros of $F_{2n}(w)$ are simple:

It is easy to see from equation (1.0) that the functions $F_{2n}(w)$ obey the differential equations

$$F_{2n}^{(2n-1)}(w) = \frac{(-1)^n}{2n} w F_{2n}(w) \qquad (2.2)$$

where $F_{2n}^{(2n-1)}(w)$ denotes the (2n-1)th derivative with respect to w.

Differentiating again yields: $F_{2n}^{(2n)}(w) = \frac{(-1)^n}{2n}(F_{2n}(w) + w F_{2n}^{(1)}(w)) \qquad (2.3)$

### **Lemma 2**

If G is a function in L-P of the form G(w) = $c \cdot \prod_{k=1}^{\infty}(1 - w/w_k)$, and G(w*) ≠ 0 and $G^{(1)}(w^*) = 0$, then $G^{(2)}(w^*)$ cannot be zero.

### **Proof of Lemma 2**

Differentiating the natural logarithm of the modulus of G as above twice yields:

$$\frac{(G^{(1)})^2 - G^{(2)}G}{G^2} = \sum_{k=1}^{\infty}(w - w_k)^{-2}$$ and we see that if G(w*) ≠ 0 (in other words, w* ≠ $w_k$ for any k)

and $G^{(1)}(w^*) = 0$, then $-G^{(2)}(w^*)G(w^*) > 0$, so $G^{(2)}(w^*) \neq 0$.

QED.

Now observe that L-P is closed under differentiation.

Set G(w) = $F_{2n}^{(2n-2)}(w)$ If for any particular $F_{2n}(w)$ there were a multiple root at w = w*, then by (2.2) and (2.3) $F_{2n}^{(2n-1)}(w^*) = 0$ and $F_{2n}^{(2n)}(w^*) = 0$. Therefore, by the above lemma, G(w*) must be zero, or a contradiction would be obtained.

So $F_{2n}^{(2n-2)}(w^*)$ and $F_{2n}^{(2n-1)}(w^*)$ must both be zero. By the same reasoning, $F_{2n}^{(2n-3)}(w^*)$,

$F_{2n}^{(2n-4)}(w^*)$,……., $F_{2n}^{(2)}(w^*)$ must all be zero, and we know $F_{2n}^{(1)}(w^*)$ and $F_{2n}(w^*)$ are both zero too. In other words, all the derivatives of $F_{2n}(w)$, up to and including the (2n-1)th, evaluated at w = w*, are zero.

Successive differentiation of equation (2.3) shows that all the derivatives of $F_{2n}(w)$ above the (2n-1)th, evaluated at w = w*, are zero as well.

So the Taylor expansion for $F_{2n}(z)$ about z = w* is:

$$F_{2n}(z) = F_{2n}(w^*) + F_{2n}^{(1)}(w^*)(z-w^*) + \frac{F_{2n}^{(2)}(w^*)}{2!}(z-w^*)^2 + \frac{F_{2n}^{(3)}(w^*)}{3!}(z-w^*)^3 + \ldots$$

and as all the coefficients $F_{2n}^{(r)}(w^*)$ are zero, we conclude $F_{2n}(z)$ is identically zero in some disk around z = w*. This establishes a contradiction, as the zeros of holomorphic functions are isolated.

Therefore there cannot be a multiple zero of $F_{2n}(w)$, in other words all the zeros are simple.

QED.

Now let us consider a particular root on the w-axis, say w = w*, for a particular value of n. Therefore $F_{2n}(w^*) = 0$ and by the above $F_{2n}^{(1)}(w^*)$ (in other words, $R_w$ evaluated at w = w*) $\neq 0$, as all the roots are simple. We also see that at z = w*, $R_\sigma = \int_\Re t.e^{-t^{2n}}.\cos(w^*.t)dt = 0$ by symmetry, so by equation (1.9), $\frac{d\omega}{d\sigma} = 0$ on the R = 0 fieldline that crosses through z = w* when $\sigma = 0$.

The geometric significance of this result is that the R = 0 fieldline that crosses each zero on the w-axis does so with a zero gradient (ie perpendicular to the w-axis).

By way of example, here is part of the plot of the R = 0 fieldlines in the case n = 2:

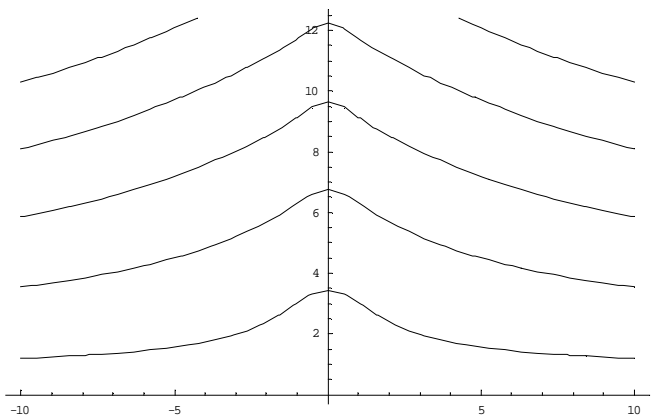

(Fig 10)

As before, the w–axis is vertical, and the $\sigma$ axis is horizontal.

Author details:   Jeremy Williams
                  12 Regent Place
                  Wimbledon
                  London SW19  8RP
                  United Kingdom

                  E-mail:   universaltutors@aol.com